\documentstyle[12pt,a4]{article}

\title{Good rotations}
\author{M. H\'enon and J-M. Petit\\
C.N.R.S., Observatoire de la C\^ote d'Azur, \\
B.P. 4229, 06304 Nice Cedex 4, France \\
{\tt henon@obs-nice.fr} and {\tt petit@obs-nice.fr}}
\date{\vskip 1truecm
Submitted to: {\it Journal of Computational Physics} \\
\medskip
May 13, 1998}

\newtheorem{theorem}{Theorem}

\newcount \mh
\ifnum \mh = 1
\input psfig.tex
\fi

\begin{document}
\maketitle

\vfill


\vspace{12pt}

\newpage

Running head: Good rotations

\vspace{36pt}
Send proofs to:
\begin{verse}
J-M. Petit\\
Observatoire de la C\^ote d'Azur \\
B.P. 4229, 06304 Nice Cedex 4, France \\
Fax: +33 4 92 00 30 33 \\
E-mail: {\tt petit@obs-nice.fr} \\
\end{verse}

\newpage

\begin{abstract}

Numerical integrations in celestial mechanics often involve the
repeated computation of a rotation with a constant angle.
A direct evaluation of these rotations yields a linear drift of
the distance to the origin.
This is due to roundoff in the representation of the sine $s$ and
cosine $c$ of the angle $\theta$.
In a computer, one generally gets $c^2 + s^2 \ne 1$, resulting in a
mapping that is slightly contracting or expanding.
In the present paper we present a method to find pairs of
representable real numbers $s$ and $c$ such that $c^2 + s^2$ is as
close to 1 as possible.
We show that this results in a drastic decrease of the systematic
error, making it negligible compared to the random error of other
operations.
We also verify that this approach gives good results in a realistic
celestial mechanics integration.

\end{abstract}

\newpage
\section{Introduction}

In some numerical computations, a rotation around a fixed axis by a
constant angle $\theta$ must be repeatedly applied.
This occurs for instance in some long-term integrations in celestial
mechanics where one must alternate between a fixed reference frame -
for integrating a Keplerian motion - and a rotating frame - to
account for some rotating perturbing potential.
A linear drift of the square distance to the axis is then generally
observed \cite{Pet98a,QT90a}, with the following properties:

\begin{itemize}

\item The rate of drift, defined as the relative change of the
  square distance to the axis per rotation, is of the order of the
  roundoff error.  For instance, if the computations are made in
  single precision, the relative change is of the order of $10 ^ {
    -7}$.

\item For a given value of $\theta$, the rate of drift is
  independent of the initial conditions.

\item The sign and the amplitude of the rate of drift seem to vary
  in quasi-random fashion with $\theta$.

\end{itemize}

These properties suggest a simple explanation for the drift.
Let us call $Z$ the rotation axis and $(X,Y)$ the plane
perpendicular to the $Z$-axis.
Then $Z$ is invariant in the rotation which is simply computed by
\begin{equation}
  \left( \begin{array}{c}
      X' \\
      Y'
    \end{array} \right)
  =
  \left( \begin{array}{cc}
      c & - s \\
      s & c
    \end{array} \right)
  \left( \begin{array}{c}
      X \\
      Y
    \end{array} \right)
                                                \label{e:simple}
\end{equation}
where ideally we should have
\begin{equation}
  c = \cos \theta, \qquad s = \sin \theta.
\end{equation}
Actually, the values of $c$ and $s$ are rounded by the computer,
and therefore $c ^ 2 + s ^ 2$ is not exactly 1. As a consequence,
the mapping (\ref{e:simple}) is slightly contracting or expanding,
in a systematic way since the same rounded values $c$ and $s$ are
used for every iteration \cite{QT90a}.

To illustrate, consider a computation in single precision, with
roundoff errors of the order of $2 ^ { -24 }$ (see below
Sect.\ref{s:roundoff}). We assume for simplicity that each step
of the computation involves one rotation. Then after $t$ steps,
the cumulative error resulting from the systematic roundoff
errors on $c$ and $s$ is $\epsilon_0 \approx 2
^ { -24 } t$. This is shown by the dotted line in
Fig.~\ref{f:errors}.

\newcommand{\legendea}{Roundoff errors as a function of the
  number of steps, in single precision.
  Dotted line: $\epsilon_0 = $ error due to the roundoff of $\cos
  \theta$ and $\sin \theta$, for an arbitrarily chosen $\theta$.
  Dashed line: $\epsilon_1 = $ error due to other roundoffs. Full
  line: $\epsilon_2 = $ errors due to the roundoff of $\cos
  \theta$ and $\sin \theta$, for a ``good rotation''
  (Eq.~(\protect\ref{e:x2+y2-3})). Dash-dot line: $\epsilon_3 = $
  same for Eq.~(\protect\ref{e:x2+y2-4}) with $k = 32$.}

\ifnum \mh = 1
\begin{figure}[htbp]
  \centerline{\psfig{figure=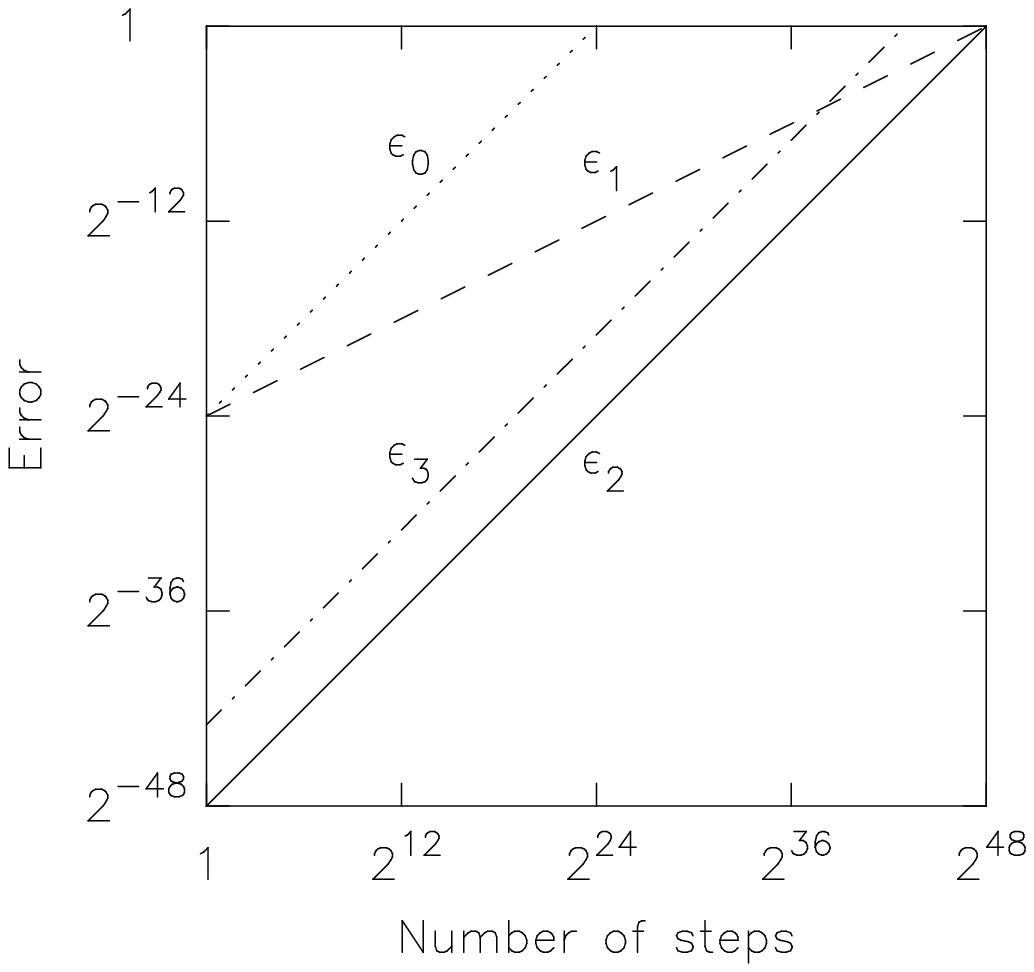,width=100mm,%
      bbllx=49mm,bblly=86mm,bburx=152mm,bbury=182mm}}
  \vspace{5mm}
  \caption{\legendea}
  \label{f:errors}
\end{figure}
\else
\begin{figure}[htbp]
  \vspace{100mm}
  \includegraphics{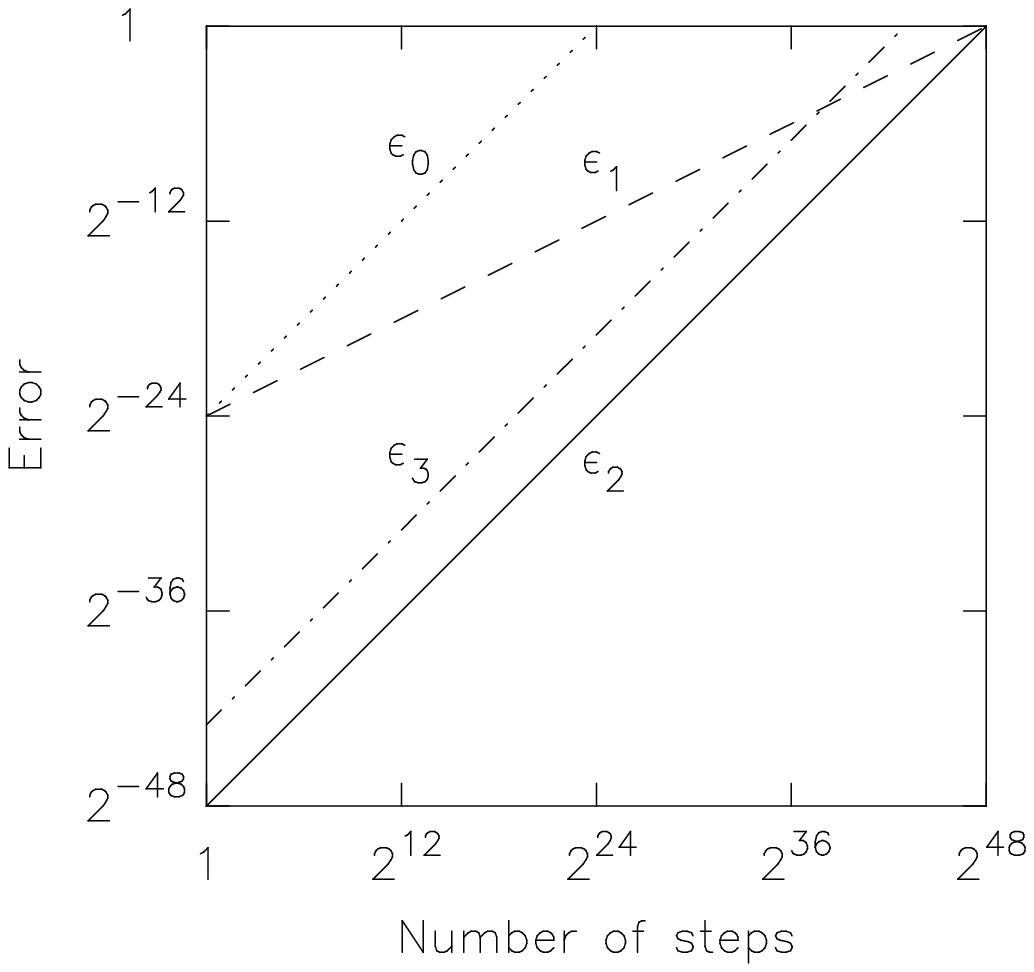}
  \begin{center}
    \caption{\legendea}
    \label{f:errors}
  \end{center}
\end{figure}
\fi

Other roundoff errors occur in the multiplications and additions
involved in (\ref{e:simple}), and in other parts of the
computation which have to be done at each step. However, these
other errors are generally quasi-random, since different values of
$X$, $Y$, and other variables are involved at each step.  A
reasonable conjecture is then that the cumulative effect has the
nature of a random walk, and that the error after $t$ steps is of
the order of $\epsilon_1 \approx 2 ^ { -24 } \sqrt { t }$. This is
represented by the dashed line in Fig.~\ref{f:errors}.

It can be seen that $\epsilon_0$ dominates. It is the cause of the
observed linear drift.

This drift can be a problem for long-term integrations. In the
case of Fig.~\ref{f:errors}, for instance, it results in a
complete breakdown of the computation after only $2 ^ { 24 }
\approx 10 ^ 7$ steps. It is therefore desirable to remove this
drift or at least to decrease its rate.

If there is some latitude in the choice of $\theta$ (for example
if it is determined by the choice of an integration step), then a
natural idea is to select a value for which the roundoff error is
very small. This is the topic of the present paper.

M. H\'enon is responsible for the mathematical basis;
J.-M. Petit, for the numerical simulations.

\section{Roundoff}
\label{s:roundoff}

A real number is usually approximated on a computer by a {\em
  representable number\/} of the form
\begin{equation}
  \sigma \times m \times 2 ^ e
\end{equation}
where $\sigma = \pm 1$ is the sign, $m$ is the mantissa and $e$ is
the exponent.

In most cases, the number is {\em normalized\/}: the exponent is
chosen in such a way that $1 / 2 \le m < 1$, i.e. the binary
representation of $m$ has the form $0.1\dots$.  The 1 in the first
position is then dropped and the next $p - 1$ binary digits are
stored. Thus, $m$ is of the form
\begin{equation}
  m = { 1 \over 2 } + { \nu \over 2 ^ p }
\end{equation}
where $\nu$ is the stored integer, which lies in the range
\begin{equation}
  0 \le \nu < 2 ^ { p - 1 }.
                                \label{e:0<=mu}
\end{equation}

Most computers today adhere to the IEEE754 standard
\cite{IEEE85a,Kah96a} and use $p = 24$ for single precision, $p
=53$ for double precision.

We consider now the binary representation $c$ of $\cos \theta$. If
$| \cos \theta | = 1$, it is exactly represented ($e = 1$, $\nu =
0$). If $1 / 2 \le | c | < 1$, the exponent is $e = 0$, and the
representable values are
\begin{equation}
  | c | = 1 / 2 + { \nu \over 2 ^ p }
\end{equation}
where $\nu$ can take all values in the range~(\ref{e:0<=mu}).  If
$1 / 4 \le | c | < 1 / 2$, the exponent is $e = - 1$,
and the representable values are
\begin{equation}
  | c | = 1 / 4 + { \nu \over 2 ^ { p + 1} }.
\end{equation}
To simplify the study, we consider only the subset of even values
of $\nu$, i.e. the values of $c$ which are multiples of $2 ^ { - p
  }$. Similarly, for $1 / 8 \le | c | < 1 / 4$, we consider only
the representable values with $\nu$ multiple of 4, and so on. In
other words, in general we consider only the representable values
of the form
\begin{equation}
  c = x 2 ^ { - p }
                                \label{e:cos-theta}
\end{equation}

where $x$ is an integer satisfying
\begin{equation}
  0 \le | x | \le 2 ^ p .
\end{equation}

Conversely, any such $x$ corresponds to a value representable on
the computer. What we have done here is simply to extract from the
cumbersome variable-size lattice of representable points a subset
of fixed size. In so doing, we eliminate some solutions of our
problem; but, as will be seen, the number of remaining solutions
is still large and should be sufficient for most applications.

The same considerations apply to $\sin \theta$, for which we
consider only the representable values of the form
\begin{equation}
  s = y 2 ^ { - p }
                                \label{e:sin-theta}
\end{equation}
where $y$ is an integer satisfying
\begin{equation}
  0 \le | y | \le 2 ^ p .
\end{equation}

\section{Some diophantine equations}

In this section we derive the equations satisfied by $x$ and $y$ for
reasonable amplitudes of the roundoff error.

\vspace{12pt}
1. We try first to find values of $\theta$ for which there is no
roundoff error, i.e. such that $\cos \theta$ and $\sin \theta$
are representable (in the restricted sense defined in the previous
Section). We are thus led to seek the solutions of the diophantine
equation
\begin{equation}
  x ^ 2 + y ^ 2 = 2 ^ { 2 p },
                                \label{e:x2+y2-1}
\end{equation}
where $p$ is given, and $x$ and $y$ are unknown integers.

Unfortunately, we have \cite{Pae90a}
\begin{theorem}
  The only solutions of~(\ref{e:x2+y2-1}) are ($x = \pm 2 ^ p$, $y
  = 0$) and ($x = 0$, $y = \pm 2 ^ p$).
\end{theorem}
We prove this recursively. If $p = 0$, the theorem is obviously
true: $x ^ 2 + y ^ 2 = 1$, so $x ^ 2 = 1$ and $y ^ 2 = 0$, or
conversely.  Assume that the theorem has been proved for $p - 1$,
with $p > 0$, and consider the value $p$. The right-hand side is
even, and $x$ and $y$ are both even or both odd. If they are both
odd, we have $x ^ 2 \bmod 4 = 1$, $y ^ 2 \bmod 4 = 1$, while $2 ^
{ 2 p } \bmod 4 = 0$: this is impossible.  If $x$ and $y$ are both
even, there is a solution $x' = x / 2$, $y' = y / 2$, for $p' = p
- 1$. According to the theorem, this solution must be of the form
($x' = \pm 2 ^ { p - 1 }$, $y' = 0$) or ($x' = 0$, $y' = \pm 2 ^ {
  p - 1 }$); from which the theorem follows.

These 4 solutions correspond to $\theta = 0$, $ \pi / 2$, $ \pi $,
$3 \pi / 2$, and are generally of no practical interest.

\vspace{12pt}
2. The next best thing which we can try is to achieve an error 1.
So we consider the diophantine equation
\begin{equation}
  x ^ 2 + y ^ 2 = 2 ^ { 2 p } - 1.
                                \label{e:x2+y2-2}
\end{equation}
As above, $p$ is given, and $x$ and $y$ are unknown integers.

But there is
\begin{theorem}
  Eq.~(\ref{e:x2+y2-2}) has no solutions for $p > 0$.
\end{theorem}
Proof: $x ^ 2 \bmod 4 = 0$ or 1, $y ^ 2 \bmod 4 = 0$ or 1, while
$2 ^ { 2 p } - 1 \bmod 4 = 3$: impossible.

\vspace{12pt}
3. So we look now for solutions of
\begin{equation}
  x ^ 2 + y ^ 2 = 2 ^ { 2 p } + 1.
                                \label{e:x2+y2-3}
\end{equation}
Fortunately, this equation always has solutions, and sometimes
many of them (see Table~\ref{t:number-of-solutions}).

The roundoff error on $c ^ 2 + s ^ 2$ is now
of the order of $2 ^ { -2 p } = 2 ^ { -48 }$ only. The cumulative
effect is $\epsilon_2 \approx 2 ^ { -48 } t$. This is represented by the
full line in Fig.~\ref{f:errors}. The situation is now inverted:
the systematic error is negligible compared to the
other errors for any realistic number of iterations. In fact both
errors become of order unity after $t = 2 ^ { 48 } \approx 3 \times
10 ^ { 14 }$ steps.

\vspace{12pt}
4. More solutions can be obtained (in order to have more choice
for the value of $\theta$), at the price of a larger roundoff
error. We look then for solutions of
\begin{equation}
  x ^ 2 + y ^ 2 = 2 ^ { 2 p } + k .
                                \label{e:x2+y2-4}
\end{equation}
This is acceptable if $k$ is not too large an integer.  The
systematic error after $t$ steps becomes $\epsilon_3 \approx 2 ^ {
  -48 } k t$. If we take for instance $k = 32$ (see
Sect.~\ref{s:single-precision}), then the error, represented by
the dash-dot line in Fig.~\ref{f:errors}, is still quite
acceptable; it becomes dominant only after $t = 2 ^ { 38 } \approx
3 \times 10 ^ { 11 }$ steps.

We give now concrete recipes for the two cases of practical
interest: single and double precision.

\section{Single precision}
\label{s:single-precision}

Because of elementary symmetries, it is clearly sufficient to
consider the range $0 \le \theta \le \pi / 4$.

We use the IEEE754 standard value, $p = 24$. Eq.~(\ref{e:x2+y2-3})
has then only 4 solutions in the range $0 \le \theta \le \pi / 4$.
Clearly this is insufficient for practical needs. So we enlarge
our search and look for solutions of~(\ref{e:x2+y2-4}), with $| k
| \le k_{\rm{max}}$.  For instance for $k_{\rm{max}} = 32$, there
are 54 solutions in the range $0 \le \theta \le \pi / 4$. These
solutions are listed in Table~{\ref{t:single}, sorted by increasing
$\theta$.

\begin{table}[p]
  \begin{center}
\leavevmode
\footnotesize
\begin{tabular}[t]{|r|r|r|r|}
  \hline
     $x$ &      $y$ &     $\theta$ &  $k$ \\
   \hline
16777216 &        0 &   0.0        &    0 \\
16777216 &        1 &    .00000006 &    1 \\
16777216 &        2 &    .00000012 &    4 \\
16777216 &        3 &    .00000018 &    9 \\
16777216 &        4 &    .00000024 &   16 \\
16777216 &        5 &    .00000030 &   25 \\
16777214 &     8192 &    .00048828 &    4 \\
16776704 &   131071 &    .00781252 &    1 \\
16761016 &   737102 &    .04394885 &    4 \\
16760654 &   745288 &    .04443725 &    4 \\
16756796 &   827505 &    .04934316 &  -15 \\
16651675 &  2048584 &    .12241060 &   25 \\
16566995 &  2647575 &    .15847021 &   -6 \\
16564945 &  2660371 &    .15924263 &   10 \\
16532686 &  2853992 &    .17094249 &    4 \\
16423326 &  3427731 &    .20575745 &  -19 \\
16389584 &  3585598 &    .21537962 &    4 \\
16332540 &  3837071 &    .23074953 &  -15 \\
16039629 &  4919886 &    .29762250 &  -19 \\
15975413 &  5124564 &    .31040869 &    9 \\
15751592 &  5776013 &    .35146882 &  -23 \\
15554306 &  6287968 &    .38417252 &    4 \\
15519136 &  6374276 &    .38972760 &   16 \\
15486659 &  6452780 &    .39479142 &   25 \\
15398649 &  6660074 &    .40821469 &   21 \\
15359229 &  6750486 &    .41409362 &  -19 \\
15263026 &  6965272 &    .42812149 &    4 \\
\hline
\end{tabular}
\nolinebreak\hspace{5mm}\nolinebreak
\begin{tabular}[t]{|r|r|r|r|}
\hline
     $x$ &      $y$ &     $\theta$ &  $k$ \\
   \hline
15259624 &  6972722 &    .42860965 &    4 \\
15045797 &  7422868 &    .45831476 &  -23 \\
14938544 &  7636418 &    .47255862 &    4 \\
14842141 &  7822137 &    .48503090 &   -6 \\
14803424 &  7895164 &    .48995756 &   16 \\
14798175 &  7904998 &    .49062199 &  -27 \\
14733952 &  8024066 &    .49868557 &    4 \\
14604001 &  8258216 &    .51464749 &    1 \\
14582860 &  8295491 &    .51720172 &   25 \\
14539039 &  8372056 &    .52245995 &    1 \\
14515158 &  8413392 &    .52530539 &  -28 \\
14362958 &  8670664 &    .54312270 &    4 \\
14188593 &  8953145 &    .56290949 &   18 \\
14067585 &  9142102 &    .57628385 &  -27 \\
13855696 &  9460162 &    .59906386 &    4 \\
13851074 &  9466928 &    .59955226 &    4 \\
13849473 &  9469270 &    .59972135 &  -27 \\
13775055 &  9577204 &    .60753567 &  -15 \\
13684899 &  9705592 &    .61688653 &    9 \\
13421774 & 10066328 &    .64350099 &    4 \\
13421771 & 10066332 &    .64350129 &    9 \\
13416856 & 10072882 &    .64398939 &    4 \\
13322259 & 10197666 &    .65332277 &  -19 \\
13012020 & 10590671 &    .68316796 &  -15 \\
12988527 & 10619470 &    .68538322 &  -27 \\
12927484 & 10693696 &    .69111140 &   16 \\
12058257 & 11665051 &    .76882501 &   -6 \\
  \hline
\end{tabular}
\normalsize
\caption{Solutions of (\protect\ref{e:x2+y2-4}) for $p = 24$, $| k
  | \le 32$.}
\label{t:single}
  \end{center}
\end{table}

This table is easily computed by scanning possible values of $y$,
which are in the range $0 \le y \le \left\lfloor \sqrt { (2 ^ { 2
      p } + k_{\rm{max}}) / 2} \right\rfloor = 11863283$; this
takes a few seconds on a workstation. (A minor technical problem
is that the terms in~(\ref{e:x2+y2-4}) are too large for the
standard integer format. This is solved by representing these
terms as double precision numbers.)

It can be seen that the values of $\theta$ cover reasonably well
the whole interval $0 \le \theta \le \pi / 4$.
If more solutions are desired, at the
expense of accepting larger roundoff errors, a larger table can
easily be built. For instance if $k_{\rm{max}} = 1000$, the
number of solutions increases to 869.

A caveat is in order here: {\em the value of $\theta$ should never
  be directly used in the computation program\/}. The values
of $\theta$ listed in Table~\ref{t:single} are not exact but
rounded; they are given here only for illustration.  Additional
unwanted roundoff would occur in computing $c$ and
$s$ from $\theta$, and the property (\ref{e:x2+y2-4})
would be destroyed in many cases.

Instead, the values of $c$ and $s$ should
receive independent names in the program, and should be computed
directly from the exact values of $x$ and $y$ listed in
Table~\ref{t:single}, using (\ref{e:cos-theta}) and
(\ref{e:sin-theta}). This computation should be done carefully,
in such a way that no roundoff occurs. In Fortran, this can be
done for instance with the instructions
\begin{verbatim}
      REAL*4  C, S
      C = 14842141. / 2. ** 24
      S = 7822137. / 2. ** 24
\end{verbatim}

\section{Double precision}
\label{s:double-precision}

For the IEEE754 standard value for double precision, $p = 53$,
Eq.~(\ref{e:x2+y2-3}) has only 8 solutions in the range $0 \le
\theta \le \pi / 4$; so we must again turn to
Eq.~(\ref{e:x2+y2-4}).

Here it is not practical to tabulate solutions
of~(\ref{e:x2+y2-4}) by scanning over $y$, as the range of
possible values of $y$ is of the order of $10 ^ { 16 }$.  Instead
we will use some classical results of the theory of numbers,
which allow a systematic generation of the solutions. We review
these results first.

\subsection{Solutions of $x ^ 2 + y ^ 2 = S$: general properties}
\label{s:x2+y2=S}

Our problem is a particular case of a more general problem: find
the solutions of the diophantine equation
\begin{equation}
  x ^ 2 + y ^ 2 = S .
                                \label{e:x2+y2-5}
\end{equation}
$S$ is a given positive integer (we disregard the trivial case $S
= 0$). This is a classical problem with a long
history \cite[Chap.~VI]{Dic20a}.

It will be convenient here to revert to a consideration of the
whole $(x, y)$ plane. We call {\em solution} a pair of integers
$x$ and $y$ satisfying~(\ref{e:x2+y2-5}).  It will also be
convenient to consider the $(x, y)$ plane as the complex plane
and to introduce the complex number
\begin{equation}
  z = x + i y = \sqrt { S } e ^ { i \theta }.
\end{equation}
(\ref{e:x2+y2-5}) can then be written
\begin{equation}
  z \bar z = S .
                                \label{e:zbarz}
\end{equation}

Note that, for a given $S$, a solution can be specified simply by
the value of $\theta$.

We call $r (S)$ the number of solutions of~(\ref{e:x2+y2-5}).
From any given solution one can deduce 3 other solutions (for the
same $S$) by rotations of $ \pi / 2$, $ \pi $, $3 \pi / 2$. In
complex notation: from any solution $z$ we deduce 3 other
solutions $i z$, $i ^ 2 z$, $ i ^ 3 z$.  These 4 solutions are
always distinct. We will call this a {\em quadruplet\/} of
solutions.

Therefore the number of solutions is a multiple of 4, and we write
$r (S) = 4 h (S)$, where $h (S)$ is the number of quadruplets.
For instance: $h (1) = 1$, $h (2) = 1$, $h (3) = 0$, $h (4) = 1$,
$h (5) = 2$, \dots

The total number of solutions up to a maximum,
\begin{equation}
  \sum_{ S = 1 } ^ { S_{\rm{max}} } r (S),
\end{equation}
is the number of points with integer coordinates inside or on the
circle of radius $\sqrt { S_{\rm{max}} }$ ; it is therefore of the
order of $ \pi S_{\rm{max}}$ \cite{Shi86a}, and the total number
of quadruplets is $\sum h (S) \sim \pi S_{\rm{max}} / 4$. From
this we deduce that the average number of quadruplets for a given
$S$ is $\langle h (S) \rangle = \pi / 4$. In practice, the
solutions are unevenly distributed. For most values of $S$, there
are no solutions. The number of values $S \le S_{\rm{max}}$ for
which $h (S) > 0$ is of the order of \cite{Shi86a}
\begin{equation}
  0.76422 { S_{\rm{max}} \over \sqrt { \log S_{\rm{max}} } } .
\end{equation}
The probability that $h (S) > 0$ for a given $S$ is obtained by
differentiating that expression:
\begin{equation}
  0.76422 \left( {1 \over \sqrt { \log S } } - {1 \over 2 (\log
      S) ^ { 3 / 2 } } \right).
\end{equation}
For values of interest here, $S \approx 2 ^ {100} \approx 10 ^
{30}$, this probability is about 0.09.

For any quadruplet generated by a solution $z$, there is a {\em
  conjugate quadruplet\/} of solutions generated by the conjugate
value $\bar z$. As is easily seen, there are 3 cases:

\begin{itemize}

\item $z$ lies on one axis, i.e. $y = 0$ or $x = 0$; $\theta \bmod
  \pi / 2 = 0$. In that case the quadruplet is identical with its
  conjugate; thus $z$ generates only 4 distinct solutions.  They
  correspond to $\theta = 0$, $ \pi / 2$, $ \pi $, $3 \pi / 2$.
  $S$ is a square in that case.

\item $z$ lies on a diagonal, i.e.  $| x | = | y |$; $\theta \bmod
  \pi / 2 = \pi / 4$. In that case again the quadruplet is
  identical with its conjugate, and $z$ generates only 4 distinct
  solutions.  They correspond to $\theta = \pi / 4$, $3 \pi / 4$,
  $5 \pi / 4$, $7 \pi / 4$. $S$ is a twice a square in that case.

\item $z$ lies neither on one axis nor on a diagonal: $\theta
  \bmod  \pi  / 4 \ne 0$. In that case the quadruplet and its
  conjugate are distinct, and $z$ generates 8 distinct solutions.
  There is one of them in each of the 8 intervals $j  \pi  / 4 <
  \theta < (j + 1)  \pi  / 4$, $j = 0, 1, \dots, 7$.

\end{itemize}

\subsection{Solutions for a given $S$}

The number of solutions for a given value of $S$ can be determined
as follows \cite[p. 242]{HW79a}. First we decompose $S$ into prime
factors.  We distinguish 3 kinds of prime factors:
\begin{itemize}
\item the factor 2,
\item factors $f_j$ equal to 1 (mod 4),
\item factors $g_j$ equal to 3 (mod 4),
\end{itemize}
and we write the decomposition of $S$ as
\begin{equation}
  S = 2 ^ \alpha \times \prod_j f_j ^ { \beta_j } \times \prod_j
  g_j ^ { \gamma_j }.
                                \label{e:S-factors}
\end{equation}
We have then the following
\begin{theorem}
  If there exists an odd $\gamma_j$, then $h (S) = 0$. If all
  $\gamma_j$ are even, then
  \begin{equation}
    h (S) = \prod_j (\beta_j + 1).
                                \label{e:h(S)}
  \end{equation}
\end{theorem}
Note that in the second case, $h (S)$ is the number of divisors of
$\prod_j f_j ^ { \beta_j }$, i.e. the number of divisors of $S$
made up of $f_j$ factors only.

We determine now the solutions themselves. 

\vspace{12pt} 1. We consider first the simple case where only one
factor $f_j$ is present, and its exponent is $\beta_j = 1$; there
are no factors 2 or $g_j$. $S$ is then a prime number equal to 1
(mod 4).  According to the above theorem, in that case $h (S) = 2$
\cite[pp. 219 and 241]{HW79a}: there are two quadruplets of
solutions.

The solutions do not lie either on an axis or on the first
diagonal ($S$ is not a square, nor twice a square). It follows
that the two quadruplets are mutually conjugate. 

We call $z_j = x_j + i y_j$ the solution with $0 < \theta < \pi /
4$ ($0 < y < x$). An algorithm exists to compute that solution for
any $f_j$ \cite{Knu81a}. The solutions for the first few factors
$f_j$ are given in Table~\ref{t:S=fj}. The two quadruplets are
generated by $z_j$ and $\bar z_j$.

\begin{table}[htbp]
  \begin{center}
    \leavevmode
\begin{tabular}{|r|r|r|r|}
      \hline
      $j$ & $f_j$ & $x_j$ & $y_j$ \\
      \hline
      1 &  5 & 2 & 1 \\
      2 & 13 & 3 & 2 \\
      3 & 17 & 4 & 1 \\
      4 & 29 & 5 & 2 \\
      5 & 37 & 6 & 1 \\
      6 & 41 & 5 & 4 \\
      7 & 53 & 7 & 2 \\
      8 & 61 & 6 & 5 \\
      9 & 73 & 8 & 3 \\
      \hline
    \end{tabular}
    \caption{Solutions in the case $S = f_j$.}
    \label{t:S=fj}
  \end{center}
\end{table}

\vspace{12pt}
2. We consider next the case where only a factor $f_j$ is
present, but with an arbitrary exponent $\beta_j$. All quadruplets
are then given by
\begin{equation}
  z = z_j ^ { \lambda_j } \bar z_j ^ { \, \beta_j - \lambda_j }
\end{equation}
where $\lambda_j$ can take the values 0, 1, \dots, $\beta_j$, and
$z_j$ is read from Table~\ref{t:S=fj}. This produces the required
number of quadruplets $h (S) = \beta_j + 1$.

Example: $S = 625 = 5 ^ 4$. Then $z_1 = 2 + i$, and the
solutions for $z$ are: $\bar z_1 ^ 4 = -7 - 24 i$, $z_1 \bar z_1 ^
3 = 15 - 20 i$, $z_1 ^ 2 \bar z_1 ^ 2 = 25$, $z_1 ^ 3 \bar z_1 =
15 + 20 i$, $z_1 ^ 4 = -7 + 24 i$.

\vspace{12pt} 3. We consider now the case with more than one
$f_j$, but still no factors 2 or $g_j$. All quadruplets are then
given by
\begin{equation}
  z = \prod_j z_j ^ { \lambda_j } \bar z_j ^ { \, \beta_j -
    \lambda_j }
                                \label{e:z-general}
\end{equation}
where $\lambda_j$ can take the values 0, 1, \dots, $\beta_j$. This
produces a number of quadruplets $h (S) = \prod_j (\beta_j + 1)$,
which is the required number.

Example: $S = 1025 = 5 ^ 2 \times 41$. There is: $f_1 = 5$,
$\beta_1 = 2$, $z_1 = 2 + i$, $f_2 = 41$, $\beta_2 = 1$, $z_2 = 5
+ 4 i$.  (\ref{e:z-general}) gives
\begin{equation}
  z = 
  \left( \begin{array}{c}
      \bar z_1 ^ 2 \\ z_1 \bar z_1 \\ z_1 ^ 2
    \end{array} \right)
  \left( \begin{array}{c}
      \bar z_2 \\ z_2
    \end{array} \right)
= \left( \begin{array}{c}
    3 - 4 i \\ 5 \\ 3 + 4 i
  \end{array} \right)
  \left( \begin{array}{c}
      5 - 4 i \\ 5 + 4 i
    \end{array} \right)
\end{equation}
where one factor should be chosen in each column.  This gives the
6 solutions $-1 -32 i$, $31 - 8 i$, $25 - 20 i$, $25 + 20 i$, $31
+ 8 i$, $-1 + 32 i$, corresponding to 6 distinct quadruplets. The
quadruplets are conjugate two by two; so there are only 3
fundamentally different solutions. In the interval $0 < \theta <
\pi / 4$, these solutions are, in terms of $x$ and $y$: $(32, 1)$,
$(31, 8)$, $(25, 20)$.

\vspace{12pt} 4. Finally, we consider the completely general case
where the exponents $\alpha$ and $\beta_j$ in (\ref{e:S-factors})
are arbitrary, and the $\gamma_j$ are even but otherwise
arbitrary. All quadruplets are then given by
\begin{equation}
  z = (1 + i) ^ \alpha \prod_j g_j ^ { \gamma_j / 2 } \prod_j
  z_j ^ { \lambda_j } \bar z_j ^ { \, \beta_j - \lambda_j }
\end{equation}

\subsection{Solutions for $S = 2 ^ { 2 n } + 1$}

In the double precision case, comparatively large values of $| k
|$ can be accepted in~(\ref{e:x2+y2-4}); even with $| k | = 10 ^
6$, for instance, the roundoff error at each step will be of the
order of $10 ^ { -26 }$ only. Thus, many more solutions can be
generated than is needed for applications. We can therefore
restrict our attention to some subset of solutions. We will
consider values of $k$ of the form $k = 2 ^ { 2 q }$, with $q \ge
0$. Consider a solution $(x, y)$ of (\ref{e:x2+y2-4}). Then $x' =
x / 2 ^ q$, $y' = y / 2 ^ q$ verify
\begin{equation}
  x' {} ^ 2 + y' {} ^ 2 = 2 ^ { 2 n } + 1
                                    \label{e:x2+y2-6}
\end{equation}
with $n = p - q$. Thus, our choice of values of $k$ is equivalent
to considering values of $S$ of the form $S (n) = 2 ^ { 2 n } +
1$, with $n \le p$.

These values have some nice properties. In particular,

\begin{itemize}
\item All prime factors of $S (n)$ are equal to 1 (mod 4). This is
  shown as follows: a prime factor $d$ of $2 ^ { 2 n } + 1$ must
  be odd. Since $2 ^ { 2 n }$ is a square, $-1$ is a quadratic
  residue (mod $d$). It follows that $(d - 1) / 2$ is even
  \cite{Ros88a}.

\item A prime factor of $S (n)$ is also a prime factor of $S (3
  n)$, $S (5 n)$, \dots This is obvious from the identity $a ^ {
    2 j + 1} + b ^ { 2 j + 1 } = (a + b) (a ^ { 2 j } - a ^ { 2 j
    - 1 } b + a ^ { 2 j - 2 } b ^ 2 - \dots + b ^ { 2 j }$, taking
  $a = 2 ^ { 2 n }$, $b = 1$ and $j = 1, 2, \dots$

\end{itemize}

As a result, the equation $x ^ 2 + y ^ 2 = S (n) = 2 ^ { 2 n } +
1$ tends to have many solutions.
Table~\ref{t:number-of-solutions} gives the number of quadruplets
$h (S)$ for $n = 1$ to 60.  This number was computed by factoring
$S$ into prime numbers (with the help of {\tt Maple}) and using
Eq.~(\ref{e:h(S)}).

\begin{table}[htbp]
\begin{center}
\leavevmode
\begin{tabular}{|r|r|}
  \hline
  $n$ & $h (S)$ \\
  \hline
   1 & 2 \\
   2 & 2 \\
   3 & 4 \\
   4 & 2 \\
   5 & 6 \\
   6 & 4 \\
   7 & 8 \\
   8 & 2 \\
   9 & 16 \\
  10 & 4 \\
  11 & 8 \\
  12 & 8 \\
  13 & 16 \\
  14 & 4 \\
  15 & 48 \\
  \hline
\end{tabular}
\hspace{10mm}
\begin{tabular}{|r|r|}
  \hline
  $n$ & $h (S)$ \\
  \hline
  16 & 4 \\
  17 & 16 \\
  18 & 16 \\
  19 & 16 \\
  20 & 4 \\
  21 & 64 \\
  22 & 8 \\
  23 & 32 \\
  24 & 8 \\
  25 & 64 \\
  26 & 8 \\
  27 & 64 \\
  28 & 8 \\
  29 & 8 \\
  30 & 16 \\
  \hline
\end{tabular}
\hspace{10mm}
\begin{tabular}{|r|r|}
  \hline
  $n$ & $h (S)$ \\
  \hline
  31 & 32 \\
  32 & 4 \\
  33 & 64 \\
  34 & 12 \\
  35 & 96 \\
  36 & 32 \\
  37 & 32 \\
  38 & 16 \\
  39 & 768 \\
  40 & 8 \\
  41 & 32 \\
  42 & 32 \\
  43 & 32 \\
  44 & 16 \\
  45 & 1536 \\
  \hline
\end{tabular}
\hspace{10mm}
\begin{tabular}{|r|r|}
  \hline
  $n$ & $h (S)$ \\
  \hline
  46 & 4 \\
  47 & 16 \\
  48 & 8 \\
  49 & 64 \\
  50 & 64 \\
  51 & 512 \\
  52 & 4 \\
  53 & 16 \\
  54 & 64 \\
  55 & 96 \\
  56 & 32 \\
  57 & 256 \\
  58 & 8 \\
  59 & 128 \\
  60 & 64 \\
  \hline
\end{tabular}
\caption{Number of quadruplets.}
\label{t:number-of-solutions}
\end{center}
\end{table}

For machines with $p = 53$, a particularly good value is
$n = 51$, for which there are 9 prime factors:
\begin{equation}
  2 ^ {102} + 1 = 1326700741 \times 26317 \times 13669 \times 3061 \times 953 \times 409 \times 137 \times 13 \times 5.
\end{equation}
Thus the total number of quadruplets is $2 ^ 9 = 512$. They are
given by the equation
\begin{eqnarray}
  x + i y =
\left( \begin{array}{c}
2 - i \\
2 + i
\end{array} \right)
\left( \begin{array}{c}
3 - 2 i \\
3 + 2 i
\end{array} \right)
\left( \begin{array}{c}
11 - 4 i \\
11 + 4 i

\end{array} \right)
\left( \begin{array}{c}
20 - 3 i \\
20 + 3 i

\end{array} \right)
\left( \begin{array}{c}
28 - 13 i \\
28 + 13 i
\end{array} \right) \nonumber \\
\left( \begin{array}{c}
55 - 6 i \\
55 + 6 i
\end{array} \right)
\left( \begin{array}{c}
113 - 30 i \\
113 + 30 i

\end{array} \right)
\left( \begin{array}{c}
154 - 51 i \\
154 + 51 i

\end{array} \right)
\left( \begin{array}{c}
30346 - 20145 i \\
30346 + 20145 i
\end{array} \right)
                                \label{e:n=51}
\end{eqnarray}
where one factor should be chosen inside each set of
parentheses. 

The angle $\theta$ is correspondingly given by
\begin{equation}
  \theta = 
\left( \begin{array}{c}
- \theta_1 \\
+ \theta_1
\end{array} \right)
+ \left( \begin{array}{c}
- \theta_2 \\
+ \theta_2
\end{array} \right)
+ \dots
+ \left( \begin{array}{c}
- \theta_{9} \\
+ \theta_{9}
\end{array} \right)
\end{equation}
with $\theta_1 = \arctan 1 / 2$, $\theta_2 = \arctan 2 / 3$, \dots
Approximate values of the $\theta_j$ are listed in
Table~\ref{t:thetaj-51}.  Here again, we point out that these
values are given only to allow an estimate of $\theta$ for a given
combination; they should never be used in the program. Instead,
the exact values of $x$ and $y$ should be computed from
(\ref{e:n=51}) for the chosen combination, and then used to
compute $c$ and $s$ as explained in
Sect.~\ref{s:single-precision}.

\begin{table}[htbp]
  \begin{center}
    \leavevmode
\begin{tabular}{|r|r|}
      \hline
      $j$ & $\theta_j$ \\
      \hline
      1 & 0.46364761 \\
      2 & 0.58800260 \\
      3 & 0.34877100 \\
      4 & 0.14888995 \\
      5 & 0.43467022 \\
      6 & 0.10866122 \\
      7 & 0.25950046 \\
      8 & 0.31980124 \\
      9 & 0.58604567 \\
      \hline
    \end{tabular}
    \caption{Values of $\theta_j$ for $n = 51$.}
    \label{t:thetaj-51}
  \end{center}
\end{table}

The total number of solutions is 2048, out of which 256 lie in
the interval $0 < \theta <  \pi   /  4$. The values of $\theta$
cover the circle quite well: the maximal difference between two
successive values is about 0.027.

Another good value is $n = 45$; there is
\begin{equation}
  2 ^ {90} + 1 = 29247661 \times 54001 \times 1321 \times 181
  \times 109 \times 61 \times 41 \times 37 \times 13 \times 5 ^ 2
\end{equation}
and the total number of quadruplets is $2 ^ 9 \times 3 = 1536$.
This value of $n$ should be appropriate in particular for machines
with $p = 48$, like some CRAYs (C90/YMP). The quadruplets are given
by the equation
\begin{eqnarray}
  x + i y =
\left( \begin{array}{c}
3 - 4 i \\
5 \\
3 + 4 i
\end{array} \right)
\left( \begin{array}{c}
3 - 2 i \\
3 + 2 i
\end{array} \right)
\left( \begin{array}{c}
6 - i \\
6 + i
\end{array} \right)
\left( \begin{array}{c}
5 - 4 i \\
5 + 4 i 
\end{array} \right)
\left( \begin{array}{c}
6 - 5 i \\
6 + 5 i 
\end{array} \right) \nonumber \\
\left( \begin{array}{c}
10 - 3 i \\
10 + 3 i 
\end{array} \right)
\left( \begin{array}{c}
10 - 9 i \\
10 + 9 i 
\end{array} \right)
\left( \begin{array}{c}
36 - 5 i \\
36 + 5 i 
\end{array} \right)
\left( \begin{array}{c}
199 - 120 i \\
199 + 120 i 
\end{array} \right)
\left( \begin{array}{c}
5331 - 910 i \\
5331 + 910 i
\end{array} \right).
                                \label{e:n=45}
\end{eqnarray}

The angle $\theta$ is correspondingly given by
\begin{equation}
  \theta = 
\left( \begin{array}{c}
- 2 \theta_1 \\
0 \\
+ 2 \theta_1
\end{array} \right)
+ \left( \begin{array}{c}
- \theta_2 \\
+ \theta_2
\end{array} \right)
+ \dots
+ \left( \begin{array}{c}
- \theta_{10} \\
+ \theta_{10}
\end{array} \right)
\end{equation}
with $\theta_1 = \arctan 1 / 2$, $\theta_2 = \arctan 2 / 3$, \dots
Approximate values of the $\theta_j$ are listed in
Table~\ref{t:thetaj-45}. 

\begin{table}[htbp]
  \begin{center}
    \leavevmode
\begin{tabular}{|r|r|}
      \hline
      $j$ & $\theta_j$ \\
      \hline
      1 & 0.46364761 \\
      2 & 0.58800260 \\
      3 & 0.16514868 \\
      4 & 0.67474094 \\
      5 & 0.69473828 \\
      6 & 0.29145679 \\
      7 & 0.73281511 \\
      8 & 0.13800602 \\
      9 & 0.54263352 \\
      10 & 0.16907011 \\
      \hline
    \end{tabular}
    \caption{Values of $\theta_j$ for $n = 45$.}
    \label{t:thetaj-45}
  \end{center}
\end{table}

The total number of solutions is 6144, out of which 768 lie in
the interval $0 < \theta <  \pi   /  4$. The values of $\theta$
cover the circle quite well: the maximal difference between two
successive values is about 0.005.

Note that some computers, like the VAXs, use a double precision
representation with $p = 56$.
However, the values of $h(S)$ for $n$ from 52 to 56 are rather small
and the values of $\theta$ cover the circle rather sparsely.

Incidentally, the mathematical approach used in the present
section could also be used in the case of single precision
(Sect.~\ref{s:single-precision}). But in that case a simple
scanning method is more convenient.

\section{Numerical simulations}
\label{s:numerical-simulations}

We now present numerical verifications of these results.  All
rotations will be performed using the traditional
mapping~(\ref{e:simple}). We point out, however, that there exist
other numerical implementations of rotations with good behaviour
over a large number of iterations \cite{Pet98a}.

\subsection{Simple rotation}

All computations will be made in double precision on a Silicon
Graphics Power Indigo 2 computer running a MIPS R8000 processor
which conforms to the IEEE754 standard for number representation.
We first study the effect of a large number
of iterations of the mapping~(\ref{e:simple}). We compare random
angles, some special angles found by chance to behave well, and
the ``good'' angles found in the previous section.

The quality of the computation is determined by the conservation
of the radius $R = \sqrt { X ^ 2 + Y ^ 2 }$. Let $R_0 = \sqrt {
  X_0 ^ 2 + Y_0 ^ 2 }$ be the radius of the initial point $(X_0,
Y_0)$. As we iterate the mapping, we record the absolute value of
the relative error $|R ^ 2 / R_0 ^ 2 - 1|$.
For each rotation angle $\theta$, the
relative error is averaged over 20 initial conditions chosen at
random. The value obtained is representative of what really
happens for all initial conditions, since the standard deviation
remains very small.

In a first series of runs, we scanned the range $0 < \theta < \pi
/ 2$ with values of the form $\theta = j / 512$, $j = 1$ to 802.
These values being representable, we can reproduce the exact same
value of $\theta$ on any computer.
Any other value of $\theta$ stored in the computer would give an
error in $c$ and $s$ of the same order of
magnitude.
However it would not be easy to kmow the exact value of $\theta$ and
reproduce the results on different computers.
Typical results are shown in Fig.~\ref{f:simple-rot}a, solid
lines. We observe a linear drift of the square radius as expected.
However, some particular angles give somewhat better results
(Fig.~\ref{f:simple-rot}a, dashed lines). For these angles, the
roundoff error on $c ^ 2 + s ^ 2$ happens to
be small and for up to $10 ^ 4$ iterations, the random errors due
to other parts of the computation are dominant. Eventually, the
linear drift emerges.  The particularly good angles presented here
correspond to $j = 126$, 248, 357, 423, and 700.

\newcommand{\legendeb}{Relative square radius errors (absolute
  value) as a function of the number of steps.
  (a): $\theta = j / 512$,
  (b): $\theta = j  \pi  / 2000$,
  (c): solutions of equation~(\protect\ref{e:x2+y2-6}) with
  $n = 45$,
  (d): solutions of equation~(\protect\ref{e:x2+y2-6}) with
  $n = 51$.}

\ifnum \mh = 1
\begin{figure}[htb]
  \centerline{\psfig{figure=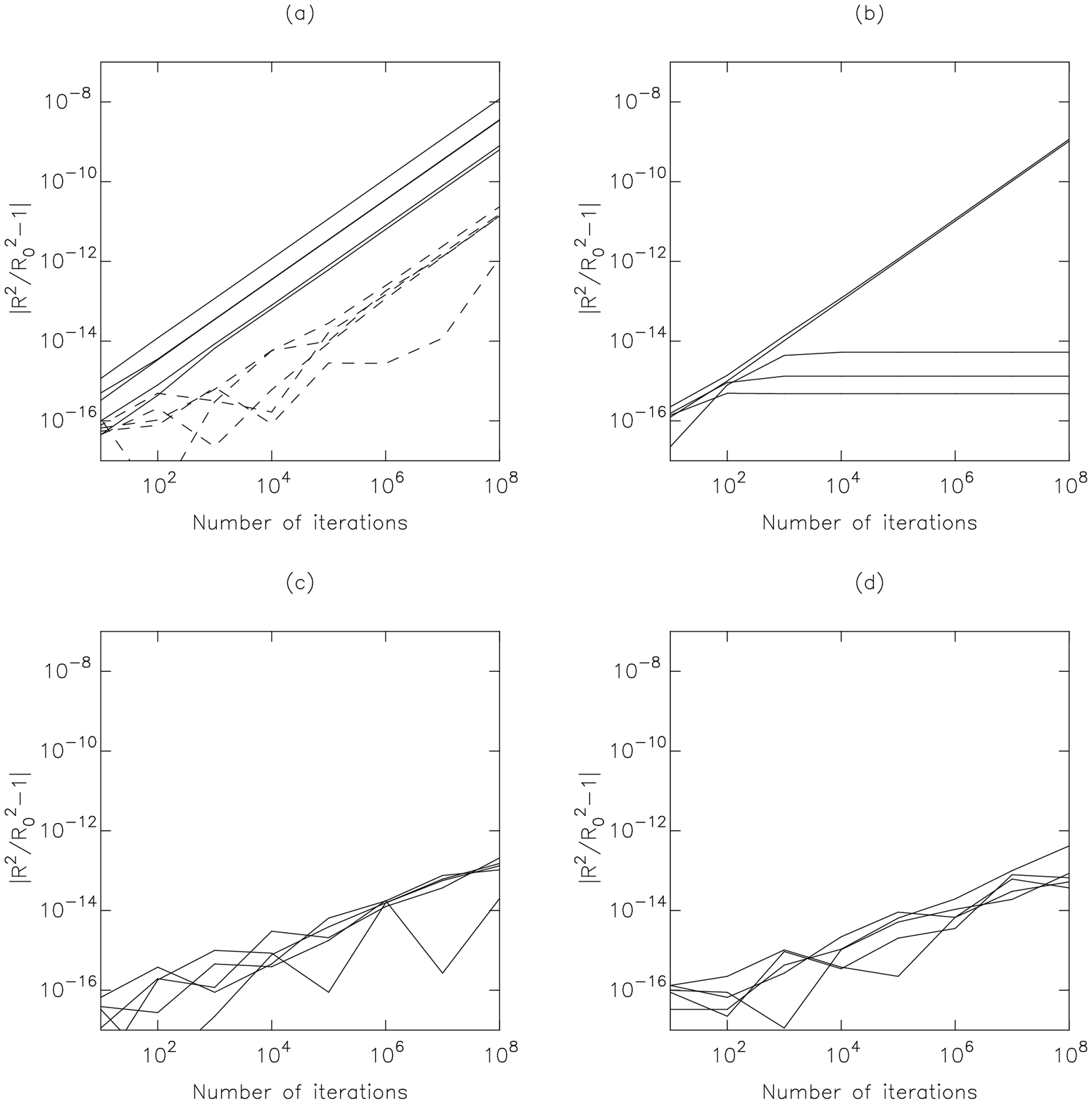,width=113mm,%
      bbllx=13mm,bblly=9mm,bburx=194mm,bbury=200mm}}
  \vspace{5mm}
  \caption{\legendeb}
  \label{f:simple-rot}
\end{figure}
\else
\begin{figure}[htbp]
  \vspace{130mm}
  \includegraphics{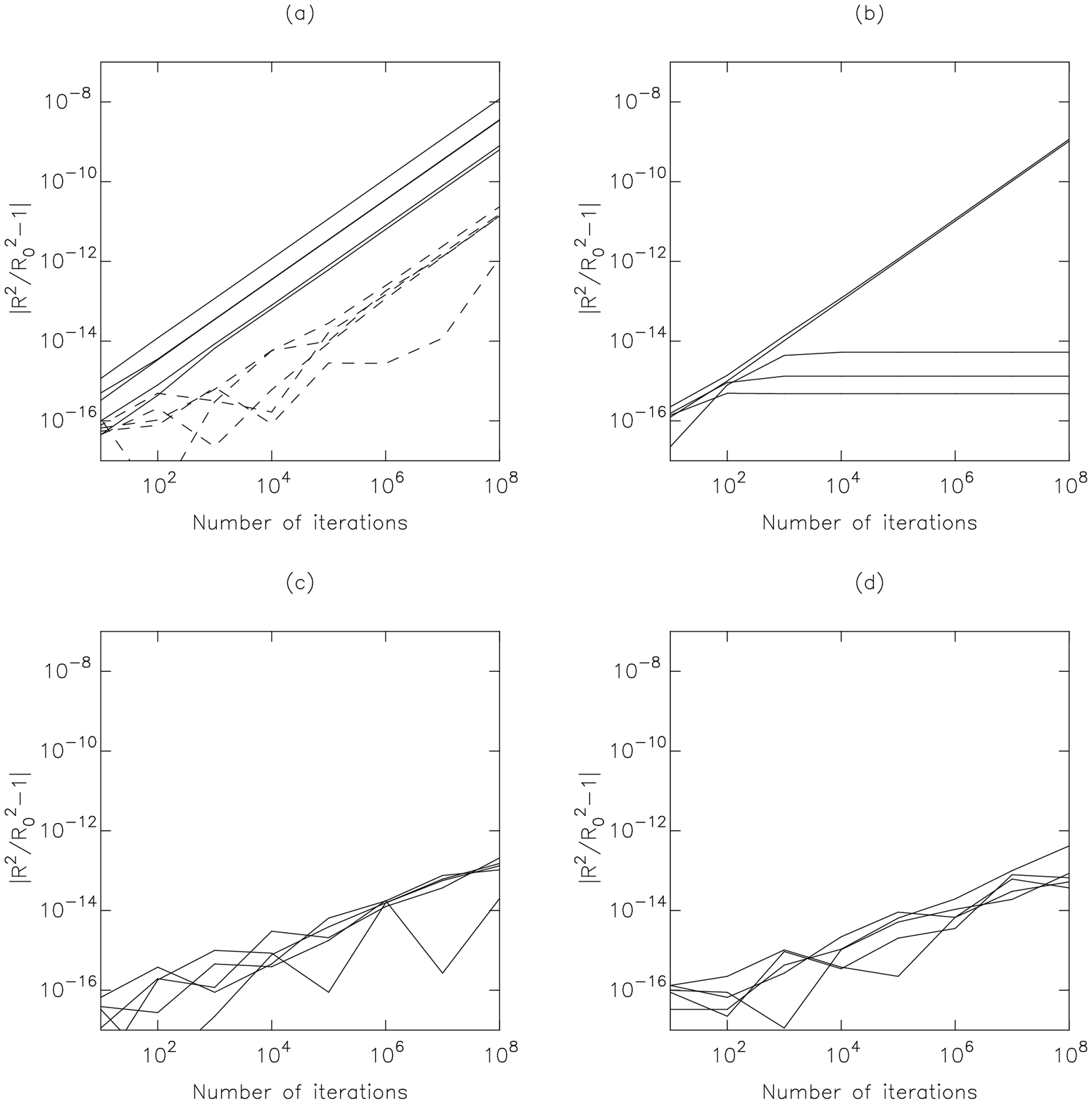}
  \begin{center}
    \caption{\legendeb}
    \label{f:simple-rot}
  \end{center}
\end{figure}
\fi

In a second series of tests, we used values of the form $\theta =
j \pi / 2000$, $j = 1$ to 999. For most values of $j$, the results
are similar to those of the previous case (Fig.~\ref{f:simple-rot}b
upper curves). However, we found 3
values of $\theta$ ($j = 40$, 250, and 450) with a peculiar
behaviour, as shown in Fig.~\ref{f:simple-rot}b (lower
curves). For $10 ^ 2$ to $10
^ 3$ iterations, the square radius drifts linearly, and then it
seems to lock on a particular value.

Inspection of the numerical results shows that a
periodic cycle of the mapping is reached. This is made
possible by the low-order commensurability of $\theta$ with $2 \pi
$. Indeed $\theta = 2 \pi / 100$, $2 \pi / 16$, and $9 \times 2
\pi / 80$ respectively and the observed periods are 100, 16, and
80.

Finally we tested the ``good'' values of $\theta$ defined by the
solutions of Eq.~(\ref{e:x2+y2-6}), with $n = 45$ and 51 (see
Table~\ref{t:test-rotation}). These solutions are integers $\le 2
^ n$. Hence, they are representable as double precision float
numbers on machines with 53 bit mantissa. Similarly $2 ^ { 45 }$
and $2 ^ { 51 }$ are representable, being powers of 2. So
assigning the solutions of Eq.~(\ref{e:x2+y2-6}) to variables and
dividing them by $2 ^ { 45 }$ or $2 ^ { 51 }$ give the desired
representable number. We next iterate the mapping. In both cases,
the random drift due to the other parts of the computation dominates
over the linear drift for longer than the $10 ^ 8$ iterations
performed here, and the overall error remains very small
(Fig.~\ref{f:simple-rot}c for $n = 45$ and
Fig.~\ref{f:simple-rot}d for $n = 51$).

\begin{table}[p]
  \begin{center}
\leavevmode
\begin{tabular}[t]{|r|r|r|}
  \hline
     \multicolumn{3}{|c|}{$n$ = 45} \\
  \hline
       $x$     &        $y$     &   $\theta$ \\
  \hline
35004143579815 &  3556679846300 & 0.10125988 \\
34476730568729 &  7021046116972 & 0.20089880 \\
33597753939071 & 10446576929072 & 0.30145465 \\
30876883071208 & 16868850912031 & 0.50001828 \\
26872087044097 & 22711912671104 & 0.70169266 \\
  \hline
     \multicolumn{3}{|c|}{$n$ = 51} \\
  \hline
        $x$      &         $y$      &   $\theta$ \\
  \hline
2240341265158844 &  226877536436263 & 0.10092511 \\
2201219968984456 &  474587240722913 & 0.21235141 \\
2150106539295032 &  669062232227809 & 0.30167850 \\
1963938109574759 & 1101612228814132 & 0.51118844 \\
1721715036961844 & 1451309661103513 & 0.70038350 \\
  \hline
\end{tabular}
\caption{Values of $x$ and $y$, and corresponding $\theta$, used in
  the numerical tests of the rotation.}
\label{t:test-rotation}
  \end{center}
\end{table}
\subsection{Integration in a rotating frame}

We came to consider this problem through the numerical study of
the long term dynamics of Dactyl, Ida's satellite
\cite{Petal97a,Pet98a}. This required the use of a symplectic
integrator in a rotating frame, thus involving a rotation. So we
want to check the effect of combining the rotation with the
iteration of the symplectic integrator of order 2 (SI2).

The implementation of SI2 we use is the generalized leap-frog
described by Yoshida \cite{Yos90a}. We write the Hamiltonian in the
form
\begin{equation}
  { \cal H } = { \cal H }_1 (L, G, H) + { \cal H }_2 (X, Y, Z) .
                                             \label{e:hamiltonian}
\end{equation}
Here, ${ \cal H }_1$ is the Hamiltonian of the two-body problem in
a rotating frame, with a primary which has the same mass as the
primary of the actual problem:
\begin{equation}
  { \cal H }_1 = - \frac { \mu ^ 2 } { 2 L ^ 2 } - \omega H,
                                             \label{e:hamil-1}
\end{equation}
$L$, $G$, and $H$ being the Delaunay variables, $\omega$ the
rotation speed of the rotating frame, and $\mu$ the product of the
gravitational constant and the reduced mass of the two bodies.
${ \cal H }_2$ represents the perturbation potential, namely
the difference between the real potential and the point mass
potential:
\begin{equation}
  { \cal H }_2 = U_{pert} (X, Y, Z) = U (X, Y, Z) + \frac { \mu }
  { R } .
                                             \label{e:hamil-2}
\end{equation}
To integrate from time $t$ to time $t + \tau$, we integrate ${
  \cal H }_2$ for $\tau / 2$, then ${ \cal H }_1$ for $\tau$ and
finally ${ \cal H }_2$ for $\tau / 2$ again.

The rotation occurs in the integration of ${ \cal H }_1$ because
of the term $- \omega H$. Using the $f$ and $g$ Gauss functions
\cite{Dan88a}, one can integrate the keplerian
Hamiltonian in a fixed frame, $- \mu ^ 2 / 2 L ^ 2$, over any time
interval $\tau$, directly in cartesian coordinates. We must then
rotate the position and velocity vectors by an angle $\omega \tau$
around the rotation axis.

Symplectic integrators are known to behave correctly on the long
run, i.e. they do not exhibit linear drifts in energy. But on a
short time scale, they may have quite large oscillating errors.
The amplitude of the oscillations are many orders of magnitude
larger than the linear drift over a period. We average the energy
error over a large number of iterations to see the secular error
rise above the oscillating error.

\newcommand{\legendec}{Relative energy errors (absolute value) as a
  function of the number of steps. (a): ``normal'' angle $\theta =
  0.0753$ (solid line) and ``good'' angles of approximately the same
  amplitude for $n = 45$ (dashed line) and $n = 51$ (dotted
  line). (b): same as (a), but for an normal angle of 0.00753, and
  corresponding good angles.}

\ifnum \mh = 1
  \begin{figure}[htbp]
  \centerline{\psfig{figure=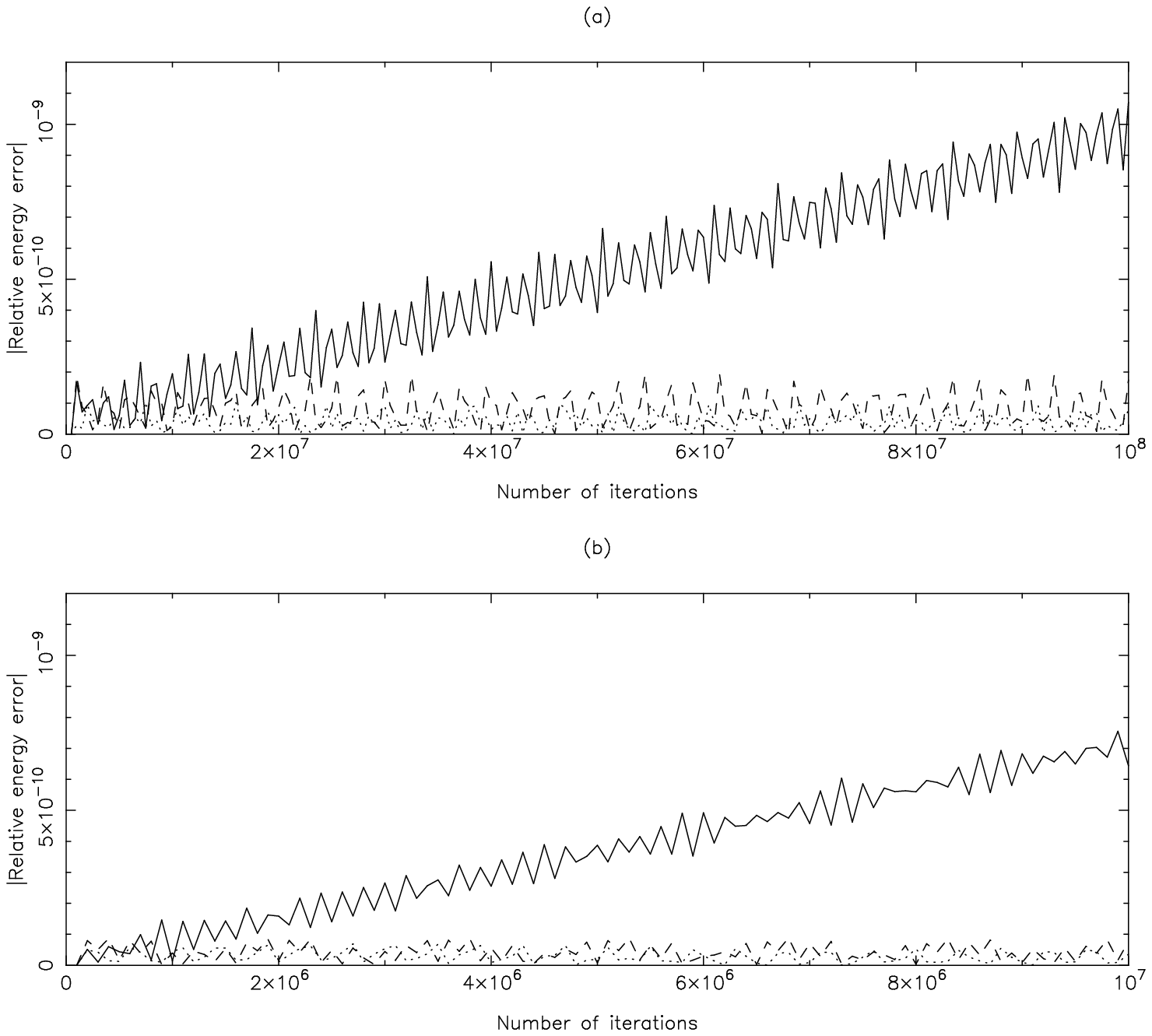,width=134mm,%
      bbllx=32mm,bblly=8mm,bburx=180mm,bbury=143mm}}
  \vspace{5mm}
  \caption{\legendec}
  \label{f:SI2}
\end{figure}
\else
\begin{figure}[htbp]
  \vspace{120mm}
  \includegraphics{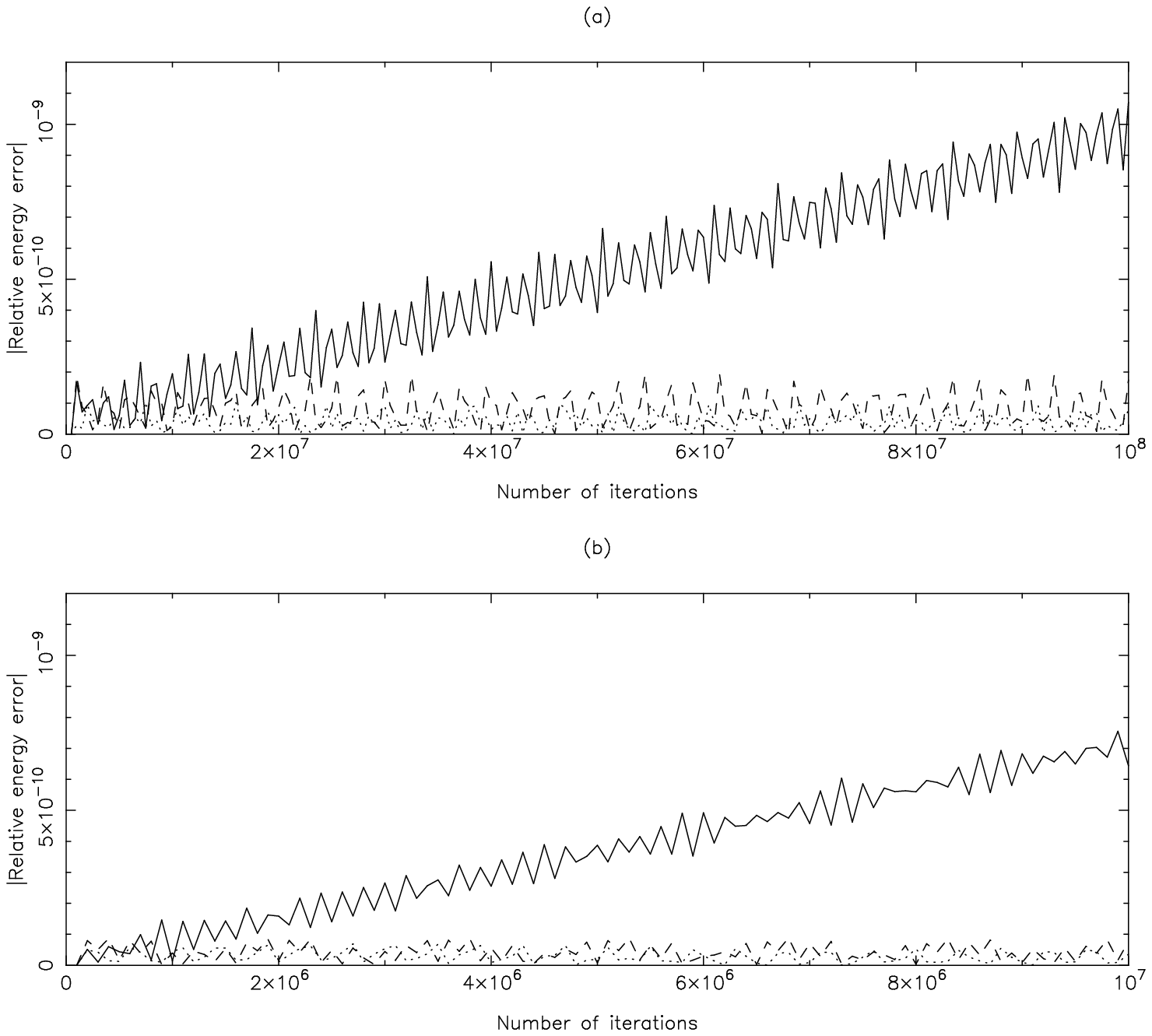}
  \begin{center}
    \caption{\legendec}
    \label{f:SI2}
  \end{center}
\end{figure}
\fi

In Fig.~\ref{f:SI2}, we present the evolution of the absolute value
of the relative energy error for two different sets of angles.
For each set, we
either choose the time step and derive the rotation angle from the
rotation speed, or we take a ``good'' angle of the same order of
magnitude for $n = 45$ or $n = 51$. Fig.~\ref{f:SI2}a shows the
error over $10 ^ 8$ iterations for an angle $\theta = 0.0753$
(solid line), $\theta \simeq 0.07511325$ (dashed line), and $\theta
\simeq 0.07194054$ (dotted line) (see
table~\ref{t:test-symplectic}, top lines). The energy error was 
averaged over $5 \times 10 ^ 5$ iterations for each data point.
For Fig.~\ref{f:SI2}b, we integrated for only $10 ^ 7$ iterations
but with an angle about 10 times smaller: $\theta = 0.00753$
(solid line), $\theta \simeq 0.00772303$ (dashed line), and $\theta
\simeq 0.00781249$ (dotted line) (see
table~\ref{t:test-symplectic}, bottom lines). For this figure, the
energy error was averaged over only $10 ^ 5$ iterations.

\begin{table}[p]
  \begin{center}
\leavevmode
\begin{tabular}[t]{|r|r|r|r|}
  \hline
$n$ &         $x$      &         $y$     &   $\theta$ \\
  \hline
 45 &   35085163629799 &   2640328077268 & 0.07511325 \\
 51 & 2245975296866668 & 161856006306841 & 0.07194054 \\
\hline
 45 &   35183322803560 &    271727410975 & 0.00772303 \\
 51 & 2251731094732799 &  17591984718848 & 0.00781249 \\
  \hline
\end{tabular}
\caption{Values of $n$, $x$ and $y$, and corresponding $\theta$,
  used in the numerical tests of the symplectic integrator.}
\label{t:test-symplectic}
  \end{center}
\end{table}

Clearly the use of solutions of Eq.~(\ref{e:x2+y2-6}) gives very
good results. We do not see any linear drift in energy.  However,
this technique can be used only if we are free to choose the
integration time step, as in the case of SI2. For example,
symplectic integrators of order 4 (SI4) or 6 as described in
\cite{Yos90a} require the use of different time steps in a very
precisely given relation: integrating over $\tau$ with SI4
corresponds to using SI2 with time step $\tau / (2 - 2 ^ { 1 / 3
  })$, then $- 2 ^ { 1 / 3 } \tau / (2 - 2 ^ { 1 / 3 })$, and
finally $\tau / (2 - 2 ^ { 1 / 3 })$ again.  This cannot be
achieved with solutions of Eq.~(\ref{e:x2+y2-6}). For such cases,
a completely different implementation of the rotation can be used,
which yields good results \cite{Pet98a}.
 
\section{Acknowledgements}

We thank F. Mignard, A. Noullez and M. Blank for discussions and
comments.


\end{document}